\documentclass[12pt]{article}%
\usepackage{amssymb}
\usepackage{amsfonts}
\usepackage{graphicx, epsfig}
\usepackage{amsmath}
\usepackage{graphicx}
\usepackage{color}%
\providecommand{\U}[1]{\protect\rule{.1in}{.1in}}
\newtheorem{theorem}{Theorem}

\begin{document}

\begin{center}
\ \ \ \ \ \ \ \ {\Large High finite-sample efficiency and robustness based on
distance-constrained maximum likelihood \vspace{1cm}}

Ricardo Maronna and V\'{\i}ctor Yohai \vspace{7cm}
\end{center}

\hrule\smallskip\smallskip

\noindent Ricardo A. Maronna is Consulting Professor, Mathematics Department,
National University of La Plata, C.C. 172, La Plata 1900, Argentina (E-mail:
rmaronna@retina.ar. V\'{\i}ctor J. Yohai is Professor Emeritus, Mathematics
Department, Faculty of Exact Sciences, Ciudad Universitaria, 1428 Buenos
Aires, Argentina (E-mail: victoryohai@gmail.com). This research was partially
supported by Grants W276 from Universidad of Buenos Aires, PIP's
112-2008-01-00216 and 112-2011-01- 00339 from CONICET and PICT 2011-0397 from
ANPCYT, Argentina.\pagebreak

\begin{center}
ABSTRACT
\end{center}

Good robust estimators can be tuned to combine a high breakdown point and a
specified asymptotic efficiency at a central model. This happens in regression
with MM- and $\tau$-estimators among others. However, the finite-sample
efficiency of these estimators can be much lower than the asymptotic one. To
overcome this drawback, an approach is proposed for parametric models, which
is based on a distance between parameters. Given a robust estimator, the
proposed one is obtained by maximizing the likelihood under the constraint
that the distance is less than a given threshold. For the linear model with
normal errors and using the MM estimator and the distance induced by the
Kullback-Leibler divergence, simulations show that the proposed estimator
attains a finite-sample efficiency close to one, while its maximum mean
squared error is smaller than that of the MM estimator. The same approach also
shows good results in the estimation of multivariate location and scatter.


\section{Introduction \label{SecIntro}}

Since the seminal work of Huber (1964) and Hampel (1971), one of the main
concerns of the research in robust statistics has been to derive statistical
procedures that are simultaneously highly robust and highly efficient under
the assumed model. The efficiency of an estimator is usually measured by the
asymptotic efficiency, that is, by the ratio between the asymptotic variances
of the maximum likelihood estimator (henceforth MLE) and of the robust
estimator. However if the sample size $n$ is not very large, this asymptotic
efficiency may be quite different from the finite sample size one, defined as
the ratio between the mean squared errors (MSE) of the MLE and of the robust
estimator, for samples of size $n$. However, it is obvious that for practical
purposes only the finite sample size efficiency matters.

Consider for example the case of a linear model with normal errors. In this
case the MLE of the regression coefficients is the least squares estimator
(LSE). It is well known that this estimator is very sensitive to outliers, and
in particular its breakdown point is zero. To overcome this problem, several
estimators combining high asymptotic breakdown point and high efficiency have
been proposed. Yohai (1987) proposed MM-estimators, which have 50\% breakdown
point and asymptotic efficiency as close to one as desired. Yohai and Zamar
(1988) proposed $\tau-$estimates, which combine the same two properties as
MM-estimators. Gervini and Yohai (2002) proposed regression estimators which
simultaneously have 50\% breakdown point and asymptotic efficiency equal to one.

However, as will be seen in Section \ref{secSimuRegre}, when $n$ is not very
large the finite sample\ efficiency of these estimators may be much smaller
than the asymptotic one. On the other hand, a $50\%$ breakdown point does not
guarantee that the estimator is highly robust. In fact, this only guarantees
that given $\varepsilon<0.5$ there exists $K(\varepsilon)$ such that if the
data are contaminated with a fraction of outliers smaller than $\varepsilon,$
the norm of the difference between the estimator and the true value is smaller
than $K(\varepsilon).$ However $K(\varepsilon)$ may be very large, which makes
the estimator unstable under outlier contamination of size $\varepsilon.$
Bondell and Stefanski (2013) proposed a regression estimator with maximum
breakdown point and high finite-sample efficiency. However, as it will be seen
in Section \ref{secSimuRegre}, the price for this efficiency is a serious loss
of robustness.

The purpose of this paper is to present estimators which have a high finite
sample size efficiency and robustness even for small $n.$ Besides, these
estimators are highly robust using a robustness criterion better than the
breakdown point, namely, the maximum MSE for a given contamination rate
$\varepsilon.$

The procedure to define the proposed estimators is very general and may be
applied to any parametric or semiparametric model. However in this paper the
details are given only to estimate the regression coefficients in a linear
model and the multivariate location and scatter of a random vector.

To define the proposed estimators we need an initial robust estimator, not
necessarily with high finite sample efficiency. Then the estimators are
defined by maximizing the likelihood function subject to the estimate being
sufficiently close to the initial one. Doing so we can expect that the
resulting estimator will\ have the maximum possible finite sample efficiency
under the assumed model compatible with proximity to the initial robust
estimator. This proximity guarantees the robustness of the new estimator.

The formulation of our proposal is as follows. Let $D$ be a distance or
discrepancy measure between densities. As a general notation, given a family
of distributions with observation vector $\mathbf{z}$, parameter vector
$\mathbf{\theta}$ and density $f\left(  \mathbf{z,\theta}\right)  ,$ put
$d\left(  \mathbf{\theta}_{1},\mathbf{\theta}_{2}\right)  =D(f\left(
\mathbf{z,\theta}_{1}\right)  ,f\left(  \mathbf{z,\theta}_{2}\right)  ).$ Let
$\mathbf{z}_{i},$ $i=1,..,n$ be i.i.d. observations with distribution
$f\left(  \mathbf{z,\theta}\right)  ,$ and let $\widehat{\mathbf{\theta}}_{0}$
be an initial robust estimator. Call $L\left(  \mathbf{z}_{1},...,\mathbf{z}%
_{n};\mathbf{\theta}\right)  $ the likelihood function. Then our proposal is
to define an estimator $\widehat{\mathbf{\theta}}$ as%
\begin{equation}
\widehat{\mathbf{\theta}}=\arg\max_{\mathbf{\theta}}L\left(  \mathbf{z}%
_{1},...,\mathbf{z}_{n};\mathbf{\theta}\right)  ~\ \mathrm{with\ }d\left(
\widehat{\mathbf{\theta}}_{0},\mathbf{\theta}\right)  \leq\delta
\label{defGeneral}%
\end{equation}
where $\delta$ is an adequately chosen constant that may depend on $n$. We
shall call this proposal \textquotedblleft distance-constrained maximum
likelihood' (DCML for short).

Several dissimilarity measures, such as the Hellinger distance, may be
employed for this purpose. We shall employ as $D$ the Kullback-Leibler (KL)
divergence, because, as it will be seen, it yields easily manageable results.
Therefore the $d$ in (\ref{defGeneral}) will be
\[
d_{\mathrm{KL}}\left(  \mathbf{\theta}_{1},\mathbf{\theta}_{2}\right)
=\int_{-\infty}^{\infty}\log\left(  \frac{f\left(  \mathbf{z,\theta}%
_{1}\right)  }{f\left(  \mathbf{z,\theta}_{2}\right)  }\right)  f\left(
\mathbf{z,\theta}_{1}\right)  d\mathbf{z.}%
\]

In Sections \ref{SecRegre} and \ref{secMulti} we apply this procedure to the
linear model and to the estimation of multivariate location and scatter,
respectively. In Section \ref{SecReal} we apply the \ DCML estimator to real data.

\section{Regression\label{SecRegre}}

Consider the family of distributions with $\mathbf{z=}\left(  \mathbf{x,}%
y\right)  ,$ with $\mathbf{x\in}R^{p}$ and $y\in R,$ satisfying the model
$y=\mathbf{x}^{\prime}\mathbf{\beta+}\sigma u,$ where $u\sim\mathrm{N}\left(
0,1\right)  $ is independent of $\mathbf{x\in}R^{p}\mathbf{.}$ Here
$\mathbf{\theta=}\left(  \mathbf{\beta,}\sigma\right)  .$ Let $\widehat
{\mathbf{\theta}}_{0}=\left(  \widehat{\mathbf{\beta}}_{0},\widehat{\sigma
}_{0}\right)  $ be an initial robust estimator of regression and scale. We
will actually consider $\sigma$ as a nuisance parameter, and therefore we
have
\begin{equation}
d_{\mathrm{KL}}\left(  \mathbf{\beta}_{0},\mathbf{\beta}\right)  =\frac
{1}{\sigma^{2}}\left(  \mathbf{\beta}-\mathbf{\beta}_{0}\right)  ^{\prime
}\mathbf{C}\left(  \mathbf{\beta}-\mathbf{\beta}_{0}\right)
~\ \label{distaRegre}%
\end{equation}
\ \textrm{with}$\mathrm{\ }\mathbf{C=}\mathrm{E}\mathbf{xx}^{\prime}.$

Here we replace $\sigma$ with its estimator $\widehat{\sigma}_{0}.$ The
natural estimator of $\mathbf{C}$ would be $\widehat{\mathbf{C}}%
=n^{-1}\mathbf{X}^{\prime}\mathbf{X}$, where $\mathbf{X}$ is the $n\times p$
matrix with rows $\mathbf{x}_{i}^{\prime}.$ Since it is not robust, we will
employ a robust version thereof. Put for $\mathbf{\beta\in}R^{p}$
$r_{i}\left(  \mathbf{\beta}\right)  =y_{i}-\mathbf{x}^{\prime}\mathbf{\beta
,}$ the residuals from $\mathbf{\beta}$. Most \textquotedblleft
smooth\textquotedblright\ robust regression estimators, like S-estimators
(Rousseeuw and Yohai, 1984) and MM-estimators satisfy the estimating equations
of an M-estimator, which can be written as weighted normal equations, namely%
\begin{equation}
\sum_{i=1}^{n}W\left(  \frac{r_{i}(\mathbf{\beta)}}{\widehat{\sigma}_{0}%
}\right)  \mathbf{x}_{i}r_{i}\left(  \mathbf{\beta}\right)  =\mathbf{0,}
\label{WNE}%
\end{equation}
where $W$ is a \textquotedblleft weight function\textquotedblright. Then we
define, as in (Yohai et al., 1991)%
\begin{equation}
\mathbf{C}_{w}=\frac{1}{\sum_{i=1}^{n}w_{i}}\sum_{i=1}^{n}w_{i}\mathbf{x}%
_{i}\mathbf{x}_{i}, \label{defCw}%
\end{equation}
with $w_{i}=W\left(  r_{i}(\widehat{\mathbf{\beta}}_{0}\mathbf{)/}%
\widehat{\sigma}_{0}\right)  .$

It is immediate that (\ref{defGeneral}) is equivalent to minimizing
$\sum_{i=1}^{n}r_{i}\left(  \mathbf{\beta}\right)  ^{2}$ subject to
$d_{\mathrm{KL}}\left(  \widehat{\mathbf{\beta}}_{0},\mathbf{\beta}\right)
\leq\delta.$ Call $\widehat{\mathbf{\beta}}_{\mathrm{LS}}$ the LSE. Put for a
general matrix $\mathbf{V:}$
\[
d_{\mathbf{V}}=\frac{1}{\widehat{\sigma}_{0}^{2}}\left(  \widehat
{\mathbf{\beta}}_{0}-\widehat{\mathbf{\beta}}_{\mathrm{LS}}\right)
\mathbf{V}\left(  \widehat{\mathbf{\beta}}_{0}-\widehat{\mathbf{\beta}%
}_{\mathrm{LS}}\right)  .
\]
Then using Lagrange multipliers, a straightforward calculation shows that
\begin{equation}
\widehat{\mathbf{\beta}}=\left\{
\begin{array}
[c]{ccc}%
\widehat{\mathbf{\beta}}_{\mathrm{LS}} & \mathrm{if} & d_{\mathbf{C}_{w}}%
\leq\delta\\
\left(  \mathbf{X}^{\prime}\mathbf{X+\lambda C}_{w}\right)  ^{-1}\left(
\mathbf{X}^{\prime}\mathbf{X}\widehat{\mathbf{\beta}}_{\mathrm{LS}%
}\mathbf{+\lambda C}_{w}\widehat{\mathbf{\beta}}_{0}\right)  & \mathrm{else}
&
\end{array}
\right.  , \label{defRegre1}%
\end{equation}
where $\lambda$ is determined from the equation $d_{\mathrm{KL}}\left(
\widehat{\mathbf{\beta}}_{0},\mathbf{\beta}\right)  =\delta$ and
$\mathbf{C}_{w}$ is defined in (\ref{defCw}). We thus see that $\widehat
{\mathbf{\beta}}$ is a linear combination of $\widehat{\mathbf{\beta}%
}_{\mathrm{0}}$ and $\widehat{\mathbf{\beta}}_{\mathrm{LS}}.$

Another approach is as follows. Define $\widehat{\mathbf{\beta}}$ as the
minimizer of $\sum_{i=1}^{n}r_{i}\left(  \mathbf{\beta}\right)  ^{2}$ subject
to $d_{\widehat{\mathbf{C}}}\leq\delta.$ In this case the solution is
explicit:
\begin{equation}
\widehat{\mathbf{\beta}}=t\widehat{\mathbf{\beta}}_{\mathrm{LS}}%
+(1-t)\widehat{\mathbf{\beta}}_{0}, \label{defRegre2}%
\end{equation}
where $t=\min\left(  1,\sqrt{\delta/d_{\widehat{\mathbf{C}}}}\right)  .$ Since
$d_{\widehat{\mathbf{C}}}$ is not robust, we now replace it with
$d_{\mathbf{C}_{w}},$ and therefore we choose%
\begin{equation}
t=\min\left(  1,\sqrt{\frac{\delta}{d_{\mathbf{C}_{w}}}}\right)  .
\label{deftRegre2}%
\end{equation}

The difference between both versions (\ref{defRegre1}) and (\ref{defRegre2})
showed to be negligible for all practical purposes.

It is easy to show that if $\widehat{\mathbf{\beta}}_{0}$ is regression- and
affine-equivariant, so is $\widehat{\mathbf{\beta}}.$

\subsection{Simulations\label{secSimuRegre}}

We now consider the model%
\begin{equation}
y_{i}=\mathbf{x}_{i}^{\prime}\mathbf{\beta+}\sigma u_{i},\ i=1,...,n,
\label{linmod}%
\end{equation}
with $\mathbf{\beta}\in R^{p}$ and $u_{i}\sim\mathrm{N}\left(  0,1\right)  $
independent of $\mathbf{x}_{i}.$ The performance of each estimator
$\widehat{\mathbf{\beta}}$ will be measured by its prediction squared error,
which is equivalent to $\left(  \widehat{\mathbf{\beta}}-\mathbf{\beta
}\right)  ^{\prime}\mathbf{C}_{x}\left(  \widehat{\mathbf{\beta}%
}-\mathbf{\beta}\right)  ,$where $\mathbf{C}_{x}\mathbf{=}\mathrm{E}%
\mathbf{xx}^{\prime}.$ Since all estimators considered are
regression-equivariant, there is no loss of generality in taking
$\mathbf{\beta=0.}$ In all cases, the distributions are normalized so that
$\mathbf{C}_{x}=\mathbf{I,}$ and therefore the criterion will be simply
$\left\Vert \widehat{\mathbf{\beta}}\right\Vert ^{2}$ where $\left\Vert
.\right\Vert $ stands for the Euclidean norm.

As initial estimator $\widehat{\mathbf{\beta}}_{0}$ we chose the MM estimator
with 85\% asymptotic efficiency and bisquare $\rho-$function:%
\begin{equation}
\rho_{\mathrm{bis}}\left(  d\right)  =1-\mathrm{I}\left(  d\leq1\right)
\left(  1-d\right)  ^{3}, \label{defbis}%
\end{equation}
\ where $\mathrm{I}\left(  .\right)  $ denotes the indicator function. The MM
estimator needs a starting regression estimator and a starting scale, which
were supplied by the S estimator $\widehat{\mathbf{\beta}}_{\mathrm{SE}}$ with
the same $\rho.$

The reason for choosing 85\% efficiency is that the maximum bias of the
resulting estimator for normal predictors is the same as that of the
regression S-estimator, as explained in Section 5.9 of (Maronna et al., 2006).

An S-estimator was also considered as an initial estimator. However, the
asymptotic efficiency of these estimators is known to be less than 33\%, and
the finite-sample efficiency is still lower. Therefore to attain acceptable
efficiencies for DCML the values $\delta$ should have to be substantially
larger than the ones we employed (given in (\ref{defidelta}) below), which
would entail a serious loss in robustness. These assertions were confirmed by
the simulations and therefore MM was the estimator of choice.

The initial scale $\widehat{\sigma}_{0}$ is a scale M estimator of the
residuals, defined as the solution of
\begin{equation}
\frac{1}{n}\sum_{i=1}^{n}\rho_{\mathrm{bis}}\left(  \frac{y_{i}-\mathbf{x}%
_{i}^{\prime}\widehat{\mathbf{\beta}}_{0}}{\mathbf{c}_{0}.\widehat{\sigma}%
_{0}}\right)  =\gamma,\label{Mscale0}%
\end{equation}
where $c_{0}=1.547$ makes  $\widehat{\sigma}_{0}$ consistent in the normal
case and $\gamma=0.5\left(  1-p/n\right)  $ 

The constant $\delta$ in (\ref{defGeneral}) is chosen as
\begin{equation}
\delta_{p,n}=0.3\frac{p}{n}. \label{defidelta}%
\end{equation}
To justify (\ref{defidelta}) note that under the model the distribution of
$nd_{\mathrm{KL}}\left(  \widehat{\mathbf{\beta}}_{0},\widehat{\mathbf{\beta}%
}_{\mathrm{LS}}\right)  $ is approximately that of $vz$ where $z\sim\chi
_{p}^{2}$ and $v$ is some constant, which implies that $\mathrm{E}%
d_{\mathrm{KL}}\left(  \widehat{\mathbf{\beta}}_{0},\widehat{\mathbf{\beta}%
}_{\mathrm{LS}}\right)  \approx vp/n.$ Therefore in order to control the
efficiency of $\widehat{\mathbf{\beta}}$ it seems reasonable to take $\delta$
of the form $cp/n$ for some $c.$ The value $c=0.3$ was arrived at after
exploratory simulations aimed at striking a balance between efficiency and robustness.

\subsubsection{Scenarios}

Since the results may depend on the distribution of the predictors, we
considered five cases, all of them including an intercept. Here each predictor
vector has the form $\mathbf{x=}\left(  1,x_{1},..,x_{p}\right)  ^{\prime
}\mathbf{,}$ where the $x_{j}$s are i.i.d. random variables with distribution
$F.$ Note that here the number of parameters is $p+1.$ In the first three
cases $F$ is standard normal, uniform in [0,1] (short-tailed) and Student with
four degrees of freedom (moderately heavy-tailed). In the other two, the
$x_{j}$s are the \emph{squares} of standard normal and uniform variables. The
Student distribution was excluded for in this case $\mathbf{C}_{x}%
\mathbf{=}\mathrm{E}\mathbf{xx}^{\prime}$ does not exist since it involves the
fourth moments of the $t_{4}$ distribution.

We took $p=5,$ 10 and 20, and $n=Kp$ with $K=5,$ 10 and 20.

For each $n$ and $p$ we first computed the finite sample efficiency. Then to
assess the estimators' robustness we contaminated the data as follows. For a
contamination rate $\varepsilon\in\left(  0,1\right)  $ let $m=[n\varepsilon]$
where $[.]$ stands for the integer part. Then for $i\leq n-m,$ $\left(
\mathbf{x}_{i},y_{i}\right)  $ were generated according to model
(\ref{linmod}), and for $i>n-m$ we put $\mathbf{x}_{i}=\left(  1,x_{0}%
,0,...,0\right)  ^{\prime}$ and $y_{i}=x_{0}K,$ where the parameter $K$ which
regulates the slope of the contamination took on a range of values in order to
determine the worst possible situations. The effect of the contamination would
be to drag the first slope towards $K.$ We took $x_{0}=5$ and $K$ ranging
between 0.5 and 2 with intervals of 0.1. We employed $\varepsilon=0.1$ and
0.2. The number of replications was $N_{\mathrm{rep}}=1000$ and $200$ for the
uncontaminated and contaminated cases respectively.

For a given scenario and estimator $\widehat{\mathbf{\beta}}$ call
$\widehat{\mathbf{\beta}}_{k},$ $k=1,...,N_{\mathrm{rep}}$ the Monte Carlo
values. As measure of performance we employed the mean squared
error:\ $\mathrm{MSE=ave}_{k}\left\{  \left\Vert \widehat{\mathbf{\beta}}%
_{k}\right\Vert ^{2}\right\}  $ where \textquotedblleft ave\textquotedblright%
\ stands for the average.

\subsubsection{Estimators}

The estimators considered were: the Least Squares estimator, the regression
S-estimator with bisquare scale (S-E), the MM estimator with bisquare loss
function and 85\% asymptotic efficiency, the Gervini-Yohai (2002) estimator
(G-Y), the Bondell-Stefanski (2013) estimator (B-S), and the proposed
estimator (DCML). Both versions (\ref{defRegre1}) and (\ref{defRegre2}) were
considered, but since the latter yielded in general slightly better results,
this is the one that is reported here. S-E, MM and G-Y were computed using the
function lmRob of the R robust package. The code for B-S was kindly supported
by the authors.

\subsubsection{Efficiency\label{SecEfiRegre}}

We deal first with the efficiencies. In order to synthesize the results, for
each combination $\left(  p,n\right)  $ we took for each estimator the minimum
efficiencies under normal errors over the five distributions, with respect to
the MLE. The results are displayed in Table \ref{tabEfi}.

\begin{center}%
\begin{table}[tbp] \centering
\bigskip%
\begin{tabular}
[c]{ccccccc}\hline
$p$ & $n$ & S-E & MM & G-Y & B-S & DCML\\\hline
\multicolumn{1}{r}{5} & \multicolumn{1}{r}{25} & \multicolumn{1}{r}{0.306} &
\multicolumn{1}{r}{0.652} & \multicolumn{1}{r}{0.657} & 0.952 &
\multicolumn{1}{r}{0.843}\\
\multicolumn{1}{r}{} & \multicolumn{1}{r}{50} & \multicolumn{1}{r}{0.270} &
\multicolumn{1}{r}{0.773} & \multicolumn{1}{r}{0.799} & 0.990 &
\multicolumn{1}{r}{0.944}\\
\multicolumn{1}{r}{} & \multicolumn{1}{r}{100} & \multicolumn{1}{r}{0.261} &
\multicolumn{1}{r}{0.810} & \multicolumn{1}{r}{0.860} & 0.996 &
\multicolumn{1}{r}{0.981}\\
\multicolumn{1}{r}{10} & \multicolumn{1}{r}{50} & \multicolumn{1}{r}{0.276} &
\multicolumn{1}{r}{0.686} & \multicolumn{1}{r}{0.702} & 0.986 &
\multicolumn{1}{r}{0.917}\\
\multicolumn{1}{r}{} & \multicolumn{1}{r}{100} & \multicolumn{1}{r}{0.276} &
\multicolumn{1}{r}{0.777} & \multicolumn{1}{r}{0.821} & 0.997 &
\multicolumn{1}{r}{0.977}\\
\multicolumn{1}{r}{} & \multicolumn{1}{r}{200} & \multicolumn{1}{r}{0.250} &
\multicolumn{1}{r}{0.808} & \multicolumn{1}{r}{0.893} & 0.999 &
\multicolumn{1}{r}{0.990}\\
\multicolumn{1}{r}{20} & \multicolumn{1}{r}{100} & \multicolumn{1}{r}{0.289} &
\multicolumn{1}{r}{0.699} & \multicolumn{1}{r}{0.723} & 0.996 &
\multicolumn{1}{r}{0.948}\\
\multicolumn{1}{r}{} & \multicolumn{1}{r}{200} & \multicolumn{1}{r}{0.254} &
\multicolumn{1}{r}{0.774} & \multicolumn{1}{r}{0.841} & 0.999 &
\multicolumn{1}{r}{0.984}\\
\multicolumn{1}{r}{} & \multicolumn{1}{r}{400} & \multicolumn{1}{r}{0.242} &
\multicolumn{1}{r}{0.820} & \multicolumn{1}{r}{0.913} & 0.999 &
\multicolumn{1}{r}{0.998}\\\hline
\end{tabular}
\caption{Minimum efficiencies of estimators for normal errors over all  x
distributions}\label{tabEfi}%
\end{table}%

\end{center}

We note the following:

\begin{itemize}
\item The efficiency of S-E is low, as can be expected

\item When $n/p$ is \textquotedblleft small\textquotedblright, the worst
finite-sample efficiency of MM can be much lower than its nominal asymptotic
one of 85\%. The worst cases with $n/p=5$ corresponded to normal
$\mathbf{x}_{i}$ with a quadratic term.

\item The worst efficiency of G-Y is also low for small $n/p.$

\item DCML\ outperforms both its initial estimator MM and G-Y.

\item B-S shows the highest efficiencies in all cases.
\end{itemize}

Table \ref{tabEfiSud} shows the efficiencies of the estimators with respect to
the MLE for model (\ref{linmod}) with Student errors $u_{i}$ with 3 and 5
degrees of freedom (\textquotedblright d.f.\textquotedblright).

\begin{center}%
\begin{table}[tbp] \centering
\begin{tabular}
[c]{cccccccc}\hline
df & $p$ & $n$ & S-E & MM & G-Y & B-S & DCML\\\hline
3 & \multicolumn{1}{r}{5} & \multicolumn{1}{r}{25} & 0.453 &
\multicolumn{1}{r}{0.828} & 0.799 & 0.875 & 0.893\\
& \multicolumn{1}{r}{} & \multicolumn{1}{r}{50} & 0.443 &
\multicolumn{1}{r}{0.917} & 0.859 & 0.883 & 0.912\\
& \multicolumn{1}{r}{} & \multicolumn{1}{r}{100} & 0.477 &
\multicolumn{1}{r}{0.949} & 0.870 & 0.871 & 0.900\\
& \multicolumn{1}{r}{10} & \multicolumn{1}{r}{50} & 0.400 &
\multicolumn{1}{r}{0.857} & 0.826 & 0.883 & 0.897\\
& \multicolumn{1}{r}{} & \multicolumn{1}{r}{100} & 0.418 &
\multicolumn{1}{r}{0.928} & 0.865 & 0.892 & 0.917\\
& \multicolumn{1}{r}{} & \multicolumn{1}{r}{200} & 0.447 &
\multicolumn{1}{r}{0.941} & 0.861 & 0.890 & 0.901\\
& \multicolumn{1}{r}{20} & \multicolumn{1}{r}{100} & 0.424 &
\multicolumn{1}{r}{0.880} & 0.854 & 0.904 & 0.943\\
& \multicolumn{1}{r}{} & \multicolumn{1}{r}{200} & 0.413 &
\multicolumn{1}{r}{0.934} & 0.881 & 0.886 & 0.904\\
& \multicolumn{1}{r}{} & \multicolumn{1}{r}{400} & 0.447 &
\multicolumn{1}{r}{0.946} & 0.867 & 0.863 & 0.883\\\hline
5 & 5 & 25 & 0.384 & 0.747 & 0.733 & 0.934 & 0.896\\
&  & 50 & 0.391 & 0.921 & 0.886 & 0.916 & 0.940\\
&  & 100 & 0.398 & 0.919 & 0.875 & 0.925 & 0.946\\
& 10 & 50 & 0.351 & 0.796 & 0.782 & 0.946 & 0.948\\
&  & 100 & 0.350 & 0.894 & 0.878 & 0.933 & 0.946\\
&  & 200 & 0.374 & 0.931 & 0.904 & 0.928 & 0.940\\
& 20 & 100 & 0.368 & 0.828 & 0.821 & 0.940 & 0.966\\
&  & 200 & 0.349 & 0.900 & 0.883 & 0.927 & 0.947\\
&  & 400 & 0.371 & 0.923 & 0.898 & 0.936 & 0.935\\\hline
\end{tabular}
\caption{Efficiencies of estimators for Student errors with 3 and 5 degrees
of freedom, and normal predictors}\label{tabEfiSud}%
\end{table}%

\end{center}

Here MM, G-Y, B-S and DCML\ exhibit high efficiencies, and none clearly
dominates the others.

\subsubsection{Robustness}

\begin{figure}
[tbh]
\begin{center}
\includegraphics[
height=9.1116cm,
width=12.1166cm
]%
{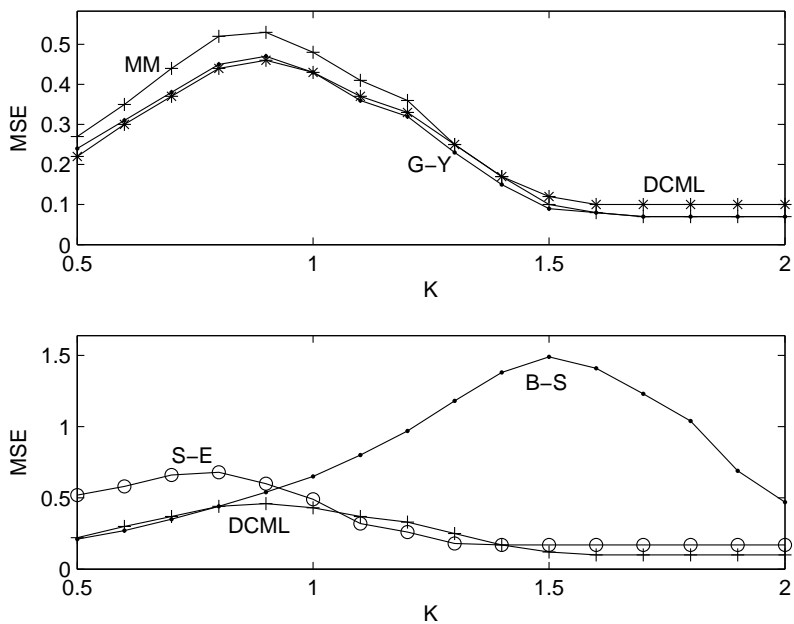}%
\caption{MSEs of regression estimators as a function of outlier size $K$ for
normal $\mathbf{x}$, $p=10,$ $n=200$ and $\varepsilon=0.1.$}%
\label{figMSEs}%
\end{center}
\end{figure}

We begin with the results of a typical case, Figure \ref{figMSEs} displays the
MSEs of the estimators for $p=10,$ $n=200,$ normal $\mathbf{x}$, and
$\varepsilon=0.1,$ for different values of the outlier size $K.$%

In the upper panel it is seen that G-Y and DCML have similar behaviors, and
that their maximum MSEs are smaller than that of MM. The lower panel shows
that the MSEs of S-E and B-S are generally larger than that of DCML, the one
of B-S being remarkably high.

Since all cases show approximately this same pattern, we display only the
maximum MSEs over $K.$ for normal $\mathbf{x}$. Table \ref{tabSimuConta} shows
the results.

\begin{center}
.\bigskip%
\begin{table}[tbp] \centering
\begin{tabular}
[c]{cccccccc}\hline
$\varepsilon$ & $p$ & $n$ & S-E & MM & G-Y & B-S & DCML\\\hline
0.1 & \multicolumn{1}{r}{5} & \multicolumn{1}{r}{25} & 1.640 &
\multicolumn{1}{r}{0.996} & 0.951 & 1.882 & 0.840\\
& \multicolumn{1}{r}{} & \multicolumn{1}{r}{50} & 1.143 &
\multicolumn{1}{r}{0.692} & 0.637 & 1.557 & 0.590\\
& \multicolumn{1}{r}{} & \multicolumn{1}{r}{100} & 0.831 &
\multicolumn{1}{r}{0.481} & 0.431 & 1.454 & 0.413\\
& \multicolumn{1}{r}{10} & \multicolumn{1}{r}{50} & 2.730 &
\multicolumn{1}{r}{1.588} & 1.514 & 2.602 & 1.268\\
& \multicolumn{1}{r}{} & \multicolumn{1}{r}{100} & 1.419 &
\multicolumn{1}{r}{0.706} & 0.644 & 1.690 & 0.597\\
& \multicolumn{1}{r}{} & \multicolumn{1}{r}{200} & 0.973 &
\multicolumn{1}{r}{0.543} & 0.475 & 1.530 & 0.463\\
& \multicolumn{1}{r}{20} & \multicolumn{1}{r}{100} & 2.058 &
\multicolumn{1}{r}{1.236} & 1.172 & 2.892 & 0.922\\
& \multicolumn{1}{r}{} & \multicolumn{1}{r}{200} & 1.212 &
\multicolumn{1}{r}{0.633} & 0.569 & 1.940 & 0.515\\
& \multicolumn{1}{r}{} & \multicolumn{1}{r}{400} & 0.850 &
\multicolumn{1}{r}{0.456} & 0.394 & 1.676 & 0.388\\\hline
0.2 & \multicolumn{1}{r}{5} & \multicolumn{1}{r}{25} & 10.51 &
\multicolumn{1}{r}{8.63} & 8.49 & 25.67 & 7.30\\
& \multicolumn{1}{r}{} & \multicolumn{1}{r}{50} & 5.24 &
\multicolumn{1}{r}{3.79} & 3.70 & 9.58 & 3.34\\
& \multicolumn{1}{r}{} & \multicolumn{1}{r}{100} & 3.17 &
\multicolumn{1}{r}{2.23} & 2.11 & 7.12 & 2.03\\
& \multicolumn{1}{r}{10} & \multicolumn{1}{r}{50} & 14.23 &
\multicolumn{1}{r}{12.00} & 11.86 & 23.37 & 9.99\\
& \multicolumn{1}{r}{} & \multicolumn{1}{r}{100} & 6.08 &
\multicolumn{1}{r}{4.10} & 3.93 & 10.84 & 3.60\\
& \multicolumn{1}{r}{} & \multicolumn{1}{r}{200} & 3.55 &
\multicolumn{1}{r}{2.47} & 2.32 & 8.70 & 2.27\\
& \multicolumn{1}{r}{20} & \multicolumn{1}{r}{100} & 6.21 &
\multicolumn{1}{r}{5.25} & 5.18 & 27.42 & 4.29\\
& \multicolumn{1}{r}{} & \multicolumn{1}{r}{200} & 3.52 &
\multicolumn{1}{r}{2.70} & 2.60 & 11.94 & 2.35\\
& \multicolumn{1}{r}{} & \multicolumn{1}{r}{400} & 2.84 &
\multicolumn{1}{r}{2.00} & 1.88 & 9.39 & 1.83\\\hline
\end{tabular}
\caption{Maximum mean squared errors of estimators with normal predictors for contaminated data}\label{tabSimuConta}%
\end{table}%

\end{center}

Some comments are in order:

\begin{itemize}
\item The MSEs of G-Y and DCML are similar, the latter being lower in most
cases. Both outperform MM, which in turn outperforms S-E.

\item The price for the high efficiency of B-S is a high contamination bias.

\item When $\varepsilon=0.2$ and $n/p=5$ all estimators have a remarkably high MSE.
\end{itemize}

As a closing comment, the joint consideration of Tables \ref{tabEfi},
\ref{tabEfiSud} and \ref{tabSimuConta} suggests that DCML shows the best
balance between efficiency and robustness.

\subsection{Asymptotic results\label{secAsyRegre}}

$\ $Assume $y=\mathbf{x}^{\prime}\mathbf{\beta+}u,$ where $u$ is independent
of $\mathbf{x}$ and has distribution $F.$ Call $\sigma_{0}$ be the limit value
of the M-scale applied to $u$ and $\mathbf{C}=\mathrm{E}(\mathbf{xx}^{\prime
}).$ It is well known that under general conditions the following expansions
hold for the MM-estimator $\widehat{\mathbf{\beta}}_{0}$ and the LS estimator
$\widehat{\mathbf{\beta}}_{\text{LS}}.$%
\[
n^{1/2}(\widehat{\mathbf{\beta}}_{0}\mathbf{-\beta)=\ }\frac{\sigma_{0}%
}{n^{1/2}\mathrm{E}\psi^{\prime}(u_{i}/\sigma_{0})}\sum_{i=1}^{n}C^{-1}%
\psi(\frac{u_{i}}{\sigma_{0}})\mathbf{x}_{i}+o\left(  \frac{1}{n^{1/2}%
}\right)  ,
\]
and
\[
n^{1/2}(\widehat{\mathbf{\beta}}_{\text{LS}}-\mathbf{\beta)=}\frac{1}{n^{1/2}%
}\sum_{i=1}^{n}C^{-1}u_{i}\mathbf{x}_{i}+o\left(  \frac{1}{n^{1/2}}\right)
\]

It then follows from the Central Limit Theorem that the joint asymptotic
distribution $J_{\mathbf{C,V}}$ of $n^{1/2}(\widehat{\mathbf{\beta}%
}_{\text{LS}}-\mathbf{\beta},\widehat{\mathbf{\beta}}-\mathbf{\beta}_{0})$ is
$J_{\mathbf{C,V}}=\mathrm{N}_{2p}$($\mathbf{0,V}\otimes\mathbf{C}^{-1})$ where
$\mathbf{V}=[V_{ij}]$ is a symmetric $2\times2$ matrix with elements
\begin{equation}
V_{11}=\mathrm{E}(u^{2}),~\ V_{12}=V_{21}=\sigma_{0}\frac{\mathrm{E}\left(
u\psi\left(  u/\sigma_{0}\right)  \right)  }{\mathrm{E}\left(  \psi^{/}\left(
u/\sigma_{0}\right)  \right)  },\ V_{22}=\sigma_{0}^{2}\frac{\mathrm{E}\left(
\psi^{2}\left(  u/\sigma_{0}\right)  \right)  }{\mathrm{E}\left(  \psi
^{/}\left(  u/\sigma_{0}\right)  \right)  } \label{VV}%
\end{equation}

Let $\left(  \mathbf{z}_{1},\mathbf{z}_{2}\right)  ^{\prime}\in R^{2p}$ be a
random vector with distribution $J_{\mathbf{C,V}}$ and define%
\begin{equation}
\mathbf{z}_{3}=t\mathbf{z}_{1}+(1-t)\mathbf{z}_{2}~\ \mathrm{with}%
t=\min\left(  1,\frac{0.3p}{(\mathbf{z}_{2}-\mathbf{z}_{1})^{\prime}%
\mathbf{C}(\mathbf{z}_{2}-\mathbf{z}_{1})}\right)  .\ \label{defz3}%
\end{equation}
Then the distribution $H_{\mathbf{C,V}}$ of $\mathbf{z}_{3}$\ is the same as
the asymptotic distribution of $n^{1/2}(\widehat{\mathbf{\beta}}%
\mathbf{-\beta}).$ Note that since $\mathbf{z}_{3}$ is a nonlinear function of
$\left(  \mathbf{z}_{1},\mathbf{z}_{2}\right)  \mathbf{,}$ $H$ is not
necessarily normal. The following Theorem will be useful determine the
distribution of $n^{1/2}\mathbf{b}^{\prime}\mathbf{(}\widehat{\mathbf{\beta}%
}\mathbf{-\beta)}$ for any $\mathbf{b\in R}^{p}$

\begin{theorem}
If $\mathbf{C=I},$ then the distribution of $v=\mathbf{d}^{\prime}%
\mathbf{z}_{3}$ is the same for any $\mathbf{d\in}R^{p}$ with $||\mathbf{d}%
||=1.$
\end{theorem}

\textbf{Proof}: Let $\mathbf{D}$ be an orthogonal matrix with first row equal
to $\mathbf{d}^{\prime}$ and let $\mathbf{v}_{j}=\mathbf{Dz}_{j},$ $1\leq
j\leq3,$ where the $\mathbf{z}_{j}$s are defined above. It is easy to check
that $(\mathbf{v}_{1},v_{2})$ has the same distribution as $(\mathbf{z}%
_{1},\mathbf{z}_{2}),$ and that $\mathbf{v}_{3}$ satisfies
\[
\mathbf{v}_{3}=t\mathbf{v}_{1}+(1-t)\mathbf{v}_{2}.
\]
Besides, we have
\[
(\mathbf{z}_{2}-\mathbf{z}_{1})^{\prime}\mathbf{C}(\mathbf{z}_{2}%
-\mathbf{z}_{1})=(\mathbf{v}_{2}-\mathbf{v}_{1})^{\prime}\mathbf{C}%
(\mathbf{v}_{2}-\mathbf{v}_{1})
\]
and therefore%
\[
t=\min\left(  1,\frac{0.3p}{(\mathbf{v}_{2}-\mathbf{v}_{1})^{\prime}%
\mathbf{C}(\mathbf{v}_{2}-\mathbf{v}_{1})}\right)  \
\]
Then $\mathbf{v}_{3}$ has the same distribution as $\mathbf{z}_{3},$ and
therefore $v_{3,1}=\mathbf{d}^{\prime}\mathbf{z}_{3}$ has the same
distribution as $z_{3,1}$ independently of $\mathbf{d.}\blacksquare$

Call $G_{\mathbf{V}}\left(  z\right)  $ the distribution function of
$v_{3,1}.$Suppose now that we want the distribution of $w=\mathbf{b}^{\prime
}\mathbf{z}_{3}$ for an arbitrary $\mathbf{C}.$ It is easy to see that
$\mathbf{z}_{3}^{\ast}=\mathbf{C}^{-1/2}z_{3}$ has distribution
$H_{\mathbf{I,V}}$ and therefore
\[
w=\mathbf{b}^{\prime}\mathbf{C}^{1/2}\mathbf{z}_{3}^{\ast}=||\mathbf{C}%
^{1/2}\mathbf{b||d}^{\prime}\mathbf{z}_{3}^{\ast}%
\]
where $||\mathbf{d}||=1.$ Then the distribution function of $w$ is
$G_{\mathbf{V}}(w/||C^{1/2}\mathbf{b||).}$

To obtain the distribution $G_{\mathbf{V}}$ we can generate a very large
sample of $\left(  \mathbf{z}_{1},\mathbf{z}_{2}\right)  $ (say of size
10$^{6}$) from $H_{\mathbf{I,V}}$ and use the transformation (\ref{defz3}) to
generate a sample of $\mathbf{z}_{3}$ with distribution $G_{\mathbf{V}}.$ In
this way we can obtain estimates of the quantiles of $G_{\mathbf{V}}$ that can
be used for asymptotic inference on any linear combination of the proposed
estimator $\widehat{\mathbf{\beta}}.$ To this end, the matrix $\mathbf{V}$ can
be estimated through (\ref{VV}), replacing $F$ by the residual empirical distribution.

This large-sample Monte Carlo can also be used to compute the asymptotic
efficiencies of $\widehat{\mathbf{\beta}}$ for different error distributions
$F.$ We compute the of $\widehat{\mathbf{\beta}}$\ with respect to the LS
estimator (eff$_{\text{LS}})$ and respect to the MM- estimator
(eff$_{\text{MM}}),$ defined by%
\[
\text{eff}_{\text{LS}}=\frac{\mathrm{E}(\mathbf{z}_{1}^{\prime}\mathbf{Cz}%
_{1})}{\mathrm{E}(\mathbf{z}_{3}^{\prime}\mathbf{Cz}_{3})},\text{
eff}_{\text{MM}}=\frac{\mathrm{E}(\mathbf{z}_{2}^{\prime}\mathbf{Cz}_{2}%
)}{\mathrm{E}(\mathbf{z}_{3}^{\prime}\mathbf{Cz}_{3})}%
\]
Since $\mathbf{z}_{1},\mathbf{z}_{2}$ and $\mathbf{z}_{3}$ are spheric when
$\mathbf{C=I}$, these efficiencies do not depend on $\mathbf{C}.$ We compute
these efficiencies when $F$ is normal, Student t with 3 and 5 degrees of
freedom, and uniform. For $p$ we chose the values $5,10$ and $20$ The results
are shown in Table \ref{tabAsymEfici}.\ Finally using the same sample we also
compute the probabilities that $\widehat{\mathbf{\beta}}$ coincides with
$\widehat{\mathbf{\beta}}_{LS}$ The results are shown in Table
\ref{tabAsymEfici}

\begin{center}%
\begin{table}[tbp] \centering
\begin{tabular}
[c]{lrrrlrrrrrrr}\hline
& \multicolumn{3}{c}{eff$_{\text{LS}}$} &  &
\multicolumn{3}{c}{eff$_{\text{MM}}$} &  &
\multicolumn{3}{c}{eff$_{\mathrm{ML-t}}$}\\\cline{1-4}\cline{6-12}%
\multicolumn{1}{r}{} & $p=5$ & $10$ & $20$ & \multicolumn{1}{r}{} & $5$ & $10
$ & $20$ &  & 5 & 10 & 20\\\cline{1-4}\cline{6-12}\cline{6-12}%
\multicolumn{1}{r}{Normal} & 0.998 & 0.9997 & 0.9999 & \multicolumn{1}{r}{} &
1.18 & 1.18 & 1.18 &  &  &  & \\
\multicolumn{1}{r}{t$_{3}$} & 1.84 & 1.84 & 1.84 & \multicolumn{1}{r}{} &
0.97 & 0.97 & 0.97 &  & 0.92 & 0.92 & 0.92\\
\multicolumn{1}{r}{t$_{5}$} & 1.19 & 1.19 & 1.19 & \multicolumn{1}{r}{} &
1.01 & 1.01 & 1.01 &  & 0.95 & 0.95 & 0.95\\
\multicolumn{1}{r}{Uniform} & 1.00 & 1.00 & 1.00 & \multicolumn{1}{r}{} &
1.07 & 1.07 & 1.07 &  &  &  & \\\hline
\end{tabular}
\caption{ Asymptotic efficiency of the proposed estimator for four error
distributions}\label{tabAsymEfici}%
\end{table}%

\end{center}

Finally\ using the same sample we also computed the probabilities that
$\widehat{\mathbf{\beta}}$ coincides with $\widehat{\mathbf{\beta}}_{LS}$ The
results are shown in Table \ref{tabProbaIgual}

\begin{center}%
\begin{table}[tbp] \centering
\begin{tabular}
[c]{rrrr}\hline
$p=$ & $5$ & $10$ & $20$\\\hline
normal & 0.85 & 0.91 & 0.96\\
t$_{3}$ & 0.02 & 0.001 & 0.00\\
t$_{5}$ & 0.14 & 0.05 & 0.01\\
uniform & 1.00 & 1.00 & 1.00\\\hline
\end{tabular}
\caption{ Probability of equality of DCML and LS estimators}\label{tabProbaIgual}%
\end{table}%

\end{center}

\subsection{Breakdown point\label{secBPRegre}}

It will be shown that for the estimators employed in this paper, the
finite-sample replacement breakdown point of the DCML estimator $\widehat
{\mathbf{\beta}}$ is that of the initial estimator $\widehat{\mathbf{\beta}%
}_{0}$.

Consider a data set $\mathbf{Z=\{z}_{i},$ $i=1,...,n\}$ with $\mathbf{z}%
_{i}=\left(  \mathbf{x}_{i},y_{i}\right)  .$ Let $m$ be such that
$\varepsilon=m/n$ is less than the breakdown point $\varepsilon^{\ast}$ of
$\widehat{\mathbf{\beta}}_{0}.$ Let $S$ (the \textquotedblleft outlier
set\textquotedblright) be any set of size $m$ contained in $\{1,...,n\}.$ Let
$\mathbf{Z}^{\ast}=\{\mathbf{z}_{i}^{\ast},~i=1,..,n\}$ where $\mathbf{z}%
_{i}^{\ast}=\mathbf{z}_{i}$ for $i\notin S$ and is arbitrary for $i\in S.$ We
have to prove that $\widehat{\mathbf{\beta}}$ is bounded as a function of
$\mathbf{Z}^{\ast}.$ The following assumptions will be needed.

A) The initial scale $\widehat{\sigma}_{0}$ is a scale M estimator of the form%
\[
\frac{1}{n}\sum_{i=1}^{n}\rho\left(  \frac{y_{i}-\mathbf{x}_{i}^{\prime
}\widehat{\mathbf{\beta}}_{0}}{\widehat{\sigma}_{0}}\right)  =\gamma,
\]
where $\rho$ is a \textquotedblleft bounded $\rho$-function\textquotedblright%
\ in the sense of (Maronna et al, 2006, p 31), i.e., $\rho\in\lbrack0,1],$
$\rho\left(  0\right)  =0,$ and $\rho\left(  t\right)  $ is a nondecreasing
function of $|t|,$ which is strictly increasing for $t>0$ such that
$\rho\left(  t\right)  <1.$

B) The breakdown point of $\widehat{\sigma}_{0}$ is $\geq\varepsilon^{\ast}.$

C) The weight function $W\left(  t\right)  $ in (\ref{WNE}) is a nondecreasing
function of $|t|$ which is \textquotedblleft matched\textquotedblright\ to
$\rho$ in the sense that $W\left(  t\right)  =0$ implies $\rho\left(
t\right)  =1.$ This is the case in the situations considered here, where
$\rho\left(  t\right)  =\rho_{\mathrm{bis}}\left(  t/\mathbf{c}_{0}.\right)  $
(see (\ref{defbis})-(\ref{Mscale0})) and $W(t)=\rho_{\mathrm{bis}}^{\prime
}\left(  t/\mathbf{c}_{1}.\right)  /t,$ where $\mathbf{c}_{1}$.$>\mathbf{c}%
_{0}$. is chosen to control the efficiency of the MM estimator.

D)\ Finally we assume%
\begin{equation}
n\left(  1-\varepsilon^{\ast}-\gamma\right)  \geq p \label{pgama}%
\end{equation}
with $\gamma$ in (\ref{Mscale0}).

Call $h$ the maximum number of $\mathbf{x}_{i}$s in a subspace. The maximal
breakdown point for $\widehat{\mathbf{\beta}}_{0}$ and $\widehat{\sigma}_{0}$
is: $\varepsilon_{\max}^{\ast}=0.5\left(  n-h-1\right)  /n.$ Here we have
$\gamma=0.5(n-p)/n\leq$ $\varepsilon_{\max}^{\ast}$ since $h\geq p-1,$ which
implies (\ref{pgama}) since $\varepsilon\leq\varepsilon_{\max}^{\ast}.$

We now proceed to the proof. Recall that $\widehat{\mathbf{\beta}}$ satisfies
\[
\frac{1}{\widehat{\sigma}_{0}^{2}}\left(  \widehat{\mathbf{\beta}}%
-\widehat{\mathbf{\beta}}_{0}\right)  ^{\prime}\mathbf{C}_{w}\left(
\widehat{\mathbf{\beta}}-\widehat{\mathbf{\beta}}_{0}\right)  \leq\delta,
\]
where $\mathbf{C}_{w}$ is defined in (\ref{defCw}). Recall that $\widehat
{\mathbf{\beta}}_{0},$ $\widehat{\sigma}_{0}$ and $\mathbf{C}_{w}$ depend on
$\mathbf{Z}^{\ast}.$ Since $\varepsilon<\varepsilon^{\ast}$ there exist
constants $a,b,c$ such that for all $S$ and $\mathbf{Z}^{\ast}:$%
\[
0<a\leq\widehat{\sigma}_{0}\leq b,~\left\Vert \widehat{\mathbf{\beta}}%
_{0}\right\Vert \leq c.
\]

Also, since $\varepsilon<\varepsilon^{\ast}$ there exists $\eta\in\left(
0,1\right)  $ such that
\begin{equation}
n\left(  1-\varepsilon-\frac{\gamma}{1-\eta}\right)  \geq p. \label{def-eta}%
\end{equation}
Let $t_{0}>0$ be such that $\rho\left(  t_{0}\right)  =1-\eta,$ and put
$w_{0}=W\left(  t_{0}\right)  .$ Then by (C) $|t|\leq t_{0}$ implies $W\left(
t\right)  \geq w_{0}>0.$ Let
\[
N=N\left(  \mathbf{Z}^{\ast}\right)  =\#\left\{  i\notin S:\rho\left(
\frac{y_{i}-\mathbf{x}_{i}^{\prime}\widehat{\mathbf{\beta}}_{0}}%
{\widehat{\sigma}_{0}}\right)  \leq1-\eta\right\}  .
\]
Then it follows from (\ref{Mscale0}) that%
\[
n\delta\geq\sum_{i\notin S}^{{}}\rho\left(  \frac{y_{i}-\mathbf{x}_{i}%
^{\prime}\widehat{\mathbf{\beta}}_{0}}{\widehat{\sigma}_{0}}\right)
\geq(n-m-N)\left(  1-\eta\right)  ,
\]
and therefore by (\ref{def-eta}), since $\varepsilon<\varepsilon^{\ast}$%
\[
N\left(  \mathbf{Z}^{\ast}\right)  \geq n-n\varepsilon-\frac{n\gamma}{1-\eta
}\geq p\forall\ \mathbf{Z}^{\ast}.
\]

Call $\mathcal{A}$ the set of all subsets of $\{1,...,n\}$ of size $h+1.$ Put%
\[
\lambda_{0}=\min_{A\in\mathcal{A}}\lambda_{\min}\left(  \sum_{i\in
A}\mathbf{x}_{i}\mathbf{x}_{i}^{\prime}\right)  ,
\]
where $\lambda_{\min}$ denotes the smallest eigenvalue of a matrix. Then
$\lambda_{0}>0$. For any vector $\mathbf{a}$ and all $\mathbf{Z}^{\ast}$ we
have%
\[
\mathbf{a}^{\prime}\mathbf{C}_{w}\mathbf{a\geq a}^{\prime}\left[
\sum_{i\notin S}W\left(  \frac{y_{i}-\mathbf{x}_{i}^{\prime}\widehat
{\mathbf{\beta}}_{0}}{\widehat{\sigma}_{0}}\right)  \mathbf{x}_{i}%
\mathbf{x}_{i}^{\prime}\right]  \mathbf{a\geq}w_{0}\lambda_{0}\left\Vert
\mathbf{a}\right\Vert ^{2},
\]
and therefore we have for all $\mathbf{Z}^{\ast}$
\[
\delta\widehat{\sigma}_{0}^{2}\geq\left(  \widehat{\mathbf{\beta}}%
-\widehat{\mathbf{\beta}}_{0}\right)  ^{\prime}\mathbf{C}_{w}\left(
\widehat{\mathbf{\beta}}-\widehat{\mathbf{\beta}}_{0}\right)  \geq
w_{0}\lambda_{0}\left\Vert \widehat{\mathbf{\beta}}-\widehat{\mathbf{\beta}%
}_{0}\right\Vert ^{2},
\]
which, in view of the boundedness of $\widehat{\mathbf{\beta}}_{0}$ and
$\widehat{\sigma}_{0},$ implies that $\widehat{\mathbf{\beta}}$ is bounded.

\section{\ Multivariate estimation\label{secMulti}}

Consider observations $\mathbf{x}_{i},$ $i=1,...,n$ with a normal $p$-variate
distribution $\mathrm{N}_{p}\left(  \mathbf{\mu,\Sigma}\right)  .$ Let
$\left(  \widehat{\mathbf{\mu}}_{0}\mathbf{,}\widehat{\mathbf{\Sigma}}%
_{0}\right)  $ be a robust estimator of multivariate location and scatter. We
shall treat $\mathbf{\mu}$ and $\mathbf{\Sigma}$ separately.

For the estimation of $\mathbf{\Sigma}$ we have, considering $\mathbf{\mu}$ as
a nuisance parameter:%
\begin{equation}
d_{\mathrm{KL}}\left(  \mathbf{\Sigma}_{0},\mathbf{\Sigma}\right)
=\log|\mathbf{\Sigma|-}\log|\mathbf{\Sigma}_{0}\mathbf{|+}\mathrm{trace}%
\left(  \mathbf{\Sigma}^{-1}\mathbf{\Sigma}_{0}\right)  -p,
\label{defDIstSIgma}%
\end{equation}
where
$\vert$%
.%
$\vert$
denotes the determinant. Our procedure amounts to
\begin{equation}
\widehat{\mathbf{\Sigma}}=\arg\min_{\mathbf{\Sigma}}\left[  n\log
|\mathbf{\Sigma|+}\sum_{i=1}^{n}\left(  \mathbf{x}_{i}-\mathbf{\mu}\right)
\mathbf{\Sigma}^{-1}\left(  \mathbf{x}_{i}-\mathbf{\mu}\right)  \right]
\label{defMultSigma}%
\end{equation}
with\textrm{ }$d_{\mathrm{KL}}\left(  \widehat{\mathbf{\Sigma}}_{0}%
,\mathbf{\Sigma}\right)  \leq\delta.$

Call $\widehat{\mathbf{\Sigma}}_{\mathrm{ML}}$ the MLE of $\mathbf{\Sigma,}$
i.e. the sample covariance matrix. Put $d_{0}=d_{\mathrm{KL}}\left(
\widehat{\mathbf{\Sigma}}_{0},\widehat{\mathbf{\Sigma}}_{\mathrm{ML}}\right)
.$ Then using Lagrange multipliers, a straightforward calculation shows that
\begin{equation}
\widehat{\mathbf{\Sigma}}=\left(  1-t\right)  \widehat{\mathbf{\Sigma}%
}_{\mathrm{ML}}+t\widehat{\mathbf{\Sigma}}_{0}, \label{defSigma_DC}%
\end{equation}
where $t=0$ if $d_{0}\leq\delta,$ and is otherwise determined from the
equation $d_{\mathrm{KL}}\left(  \widehat{\mathbf{\Sigma}}_{0},\mathbf{\Sigma
}\right)  =\delta,$ which is easily derived from (\ref{defDIstSIgma}%
)-(\ref{defSigma_DC}).

We now turn to $\mathbf{\mu.}$ We have%
\[
d_{\mathrm{KL}}\left(  \mathbf{\mu}_{0},\mathbf{\mu}\right)  =\left(
\mathbf{\mu-\mu}_{0}\right)  ^{\prime}\Sigma^{-1}\left(  \mathbf{\mu-\mu}%
_{0}\right)  .
\]
The estimator is then defined by
\begin{equation}
\sum_{i=1}^{n}\left(  \mathbf{x}_{i}\mathbf{-\mu}\right)  \mathbf{\Sigma}%
^{-1}\left(  \mathbf{x}_{i}\mathbf{-\mu}\right)  =\min\label{defMultMu}%
\end{equation}
with $d_{\mathrm{KL}}\left(  \mathbf{\mu}_{0},\mathbf{\mu}\right)  \leq
\delta.$ Let $\overline{\mathbf{x}}$ be the sample mean, and define%
\[
d_{0}=\left(  \overline{\mathbf{x}}\mathbf{-}\widehat{\mathbf{\mu}}%
_{0}\right)  ^{\prime}\widehat{\Sigma}_{0}^{-1}\left(  \mathbf{x-}%
\widehat{\mathbf{\mu}}_{0}\right)  .
\]
Then a straightforward calculation shows that
\begin{equation}
\widehat{\mathbf{\mu}}=t\overline{\mathbf{x}}\mathbf{+}\left(  1-t\right)
\widehat{\mathbf{\mu}}_{0} \label{defmu_DC}%
\end{equation}
with%
\[
t=\min\left(  1,\sqrt{\frac{\delta}{d_{0}}}\right)  .
\]

It is easy to show that if the initial estimators are affine-equivariant, so
are the resulting ones.

\textbf{Remark: }Unlike the regression and location cases, $d_{\mathrm{KL}%
}\left(  \mathbf{\Sigma}_{0},\mathbf{\Sigma}\right)  $ is not symmetric in its
arguments. Here we have chosen the form (\ref{defDIstSIgma}) because it yields
the simple intuitive result (\ref{defSigma_DC}), while the alternative order
yields a more complicated result.

\subsection{Simulations\label{secSimuMulti}}

As initial estimator we employ an S estimator (Davies, 1987)\ with bisquare
scale, computed as described at the end of page 199 of (Maronna et al, 2006).
It is  implemented  as the \ function CovSest  with the option method=
"bisquare" in the R package rrcov.

This study includes $p=2,$ 5 and 10. The reason why larger values of $p$ are
not included is the following. Rocke (1996) found out that the efficiency of S
estimators with a monotone weight function increases with $p,$ and therefore
there is little to be gained with DCML when $p$ is large.

We now define the S estimator. For $\left(  \mathbf{\mu,\Sigma}\right)  $
denote the (squared) Mahalanobis distance of $\mathbf{x}$ as
\[
d\left(  \mathbf{x,\mu,\Sigma}\right)  =\left(  \mathbf{x-\mu}\right)
^{\prime}\mathbf{\Sigma}^{-1}\left(  \mathbf{x-\mu}\right)  .
\]
Define a scale M estimator $\widehat{\sigma}=\widehat{\sigma}\left(
\mathbf{\mu,\Sigma}\right)  $ as the solution of
\[
\frac{1}{n}\sum_{i=1}^{n}\rho\left(  \frac{d\left(  \mathbf{x,\mu,\Sigma
}\right)  ^{1/2}}{\sigma}\right)  =\gamma,
\]
where $\rho$ is the bisquare $\rho$-function (\ref{defbis}), and
$\gamma=0.5\left(  1-p/n\right)  $ which ensures maximal breakdown point. The
S estimator is defined by%
\[
\left(  \widehat{\mathbf{\mu}}_{0},\widetilde{\mathbf{\Sigma}}\right)
=\arg\min\left\{  \widehat{\sigma}\left(  \mathbf{t,V}\right)  :\mathbf{t\in
}R^{p},|\mathbf{V}|=1\right\}
\]

Since $|\widetilde{\mathbf{\Sigma}}|=1,$ we have to scale $\widetilde
{\mathbf{\Sigma}}$ to make it a consistent estimator of the covariance matrix
under normality. Put $d_{i}=d\left(  \mathbf{x}_{i},\widehat{\mathbf{\mu}}%
_{0},\widetilde{\mathbf{\Sigma}}\right)  $ and call $\chi_{p}^{2}$ the
chi-squared distribution with $p$ degrees of freedom. Then define%
\[
\widehat{\mathbf{\Sigma}}_{0}=\frac{\mathrm{median}_{i}\{d_{i}\}}%
{\mathrm{median}(\chi_{p}^{2})}\widetilde{\mathbf{\Sigma}}.
\]
The constants $\delta$ in (\ref{defMultSigma}) and (\ref{defMultMu}) were
chosen as
\begin{equation}
\delta=an^{-b}p^{c}, \label{deltamulti}%
\end{equation}
with $\left(  a,b,c\right)  $ given in Table \ref{tabConstCuan}.

\begin{center}
\bigskip%
\begin{table}[tbp] \centering
$%
\begin{tabular}
[c]{cccc}\hline
& a & b & c\\\hline
\multicolumn{1}{r}{$\mathbf{\Sigma}$} & \multicolumn{1}{r}{1.02} &
\multicolumn{1}{r}{0.82} & \multicolumn{1}{r}{0.18}\\
\multicolumn{1}{r}{$\mathbf{\mu}$} & \multicolumn{1}{r}{0.55} &
\multicolumn{1}{r}{0.88} & \multicolumn{1}{r}{-0.30}\\\hline
\end{tabular}
\ $%
\caption{Constants for the approximate computation of $\delta$ }\label{tabConstCuan}%
\end{table}%

\end{center}

The motivation for this choice is as follows. It was considered as reasonable
to choose for each $\left(  p,n\right)  ,$ $\delta$ as some $\alpha-$quantile
of $d_{\mathrm{KL}}$ under the nominal model, i.e. the multivariate normal
distribution. Exploratory simulations suggested $\alpha$ between 0.4 and 0.6.
The quantiles were computed by simulation for $p$ between 2 and 10 and $n$
between $5p$ and 500. Then for each $\alpha$ the $\alpha$-quantile was fitted
by regression as a function of $n$ and $p$ of the form (\ref{deltamulti}).
Finally, after the simulation was completed, it was decided that $\alpha=0.4$
yielded the best results.

The values of $c$ indicate that when $p$ increases, the quantiles for
$\mathbf{\Sigma}$ increase very slowly, and those for $\mathbf{\mu}$ decrease.
This fact may seem counter-intuitive, but it is a consequence of the
increasing efficiency of the S estimator: when $p$ increases, the S estimator
becomes \textquotedblleft closer\textquotedblright\ to the classical one,
which makes $d_{\mathrm{KL}}$ smaller.

For each $n$ and $p$ we generate $N_{\mathrm{rep}}$ samples of size $n$ from
$\mathrm{N}_{p}\left(  \mathbf{0,I}\right)  .$ For a contamination rate
$\varepsilon,$ the first $m=[n\varepsilon]$ elements are replaced by $\left(
K,0,...,0\right)  $ where $K$ ranges between 1 and 10. For each sample three
estimators were computed: the sample mean and covariance matrix, the S
estimator, and the DCML estimator given by (\ref{defSigma_DC})-(\ref{defmu_DC}).

For each scenario, each estimator is evaluated by its \textquotedblleft
loss\textquotedblright\ defined as $\left\Vert \mathbf{\mu}\right\Vert ^{2}$
for location and as $d_{\mathrm{KL}}\left(  \mathbf{I,\Sigma}\right)
=\mathrm{trace}\left(  \mathbf{\Sigma}\right)  -\log|\mathbf{\Sigma}|$ for
scatter and the results were summarized by the respective mean losses. Table
\ref{tabEfiMulti} shows the efficiencies, defined as the ratio of the mean
losses of the classical and the robust estimator.

\begin{center}%
\begin{table}[tbp] \centering
\begin{tabular}
[c]{ccccccc}\hline
& \textbf{\ } & \multicolumn{2}{c}{\textbf{\ }$\mathbf{\Sigma}$} &  &
\multicolumn{2}{c}{$\mathbf{\mu}$}\\\hline
$p$ & $n$ & S-E & DCML &  & S-E & DCML\\\cline{1-4}\cline{6-7}%
\multicolumn{1}{r}{2} & \multicolumn{1}{r}{10} & \multicolumn{1}{r}{0.422} &
\multicolumn{1}{r}{0.627} & \multicolumn{1}{r}{} & \multicolumn{1}{r}{0.690} &
\multicolumn{1}{r}{0.889}\\
\multicolumn{1}{r}{} & \multicolumn{1}{r}{20} & \multicolumn{1}{r}{0.414} &
\multicolumn{1}{r}{0.692} & \multicolumn{1}{r}{} & \multicolumn{1}{r}{0.673} &
\multicolumn{1}{r}{0.898}\\
\multicolumn{1}{r}{} & \multicolumn{1}{r}{40} & \multicolumn{1}{r}{0.407} &
\multicolumn{1}{r}{0.762} & \multicolumn{1}{r}{} & \multicolumn{1}{r}{0.586} &
\multicolumn{1}{r}{0.867}\\
\multicolumn{1}{r}{5} & \multicolumn{1}{r}{25} & \multicolumn{1}{r}{0.772} &
\multicolumn{1}{r}{0.922} & \multicolumn{1}{r}{} & \multicolumn{1}{r}{0.893} &
\multicolumn{1}{r}{0.971}\\
\multicolumn{1}{r}{} & \multicolumn{1}{r}{50} & \multicolumn{1}{r}{0.778} &
\multicolumn{1}{r}{0.962} & \multicolumn{1}{r}{} & \multicolumn{1}{r}{0.876} &
\multicolumn{1}{r}{0.980}\\
\multicolumn{1}{r}{} & \multicolumn{1}{r}{100} & \multicolumn{1}{r}{0.777} &
\multicolumn{1}{r}{0.977} & \multicolumn{1}{r}{} & \multicolumn{1}{r}{0.855} &
\multicolumn{1}{r}{0.978}\\
\multicolumn{1}{r}{10} & \multicolumn{1}{r}{50} & \multicolumn{1}{r}{0.936} &
\multicolumn{1}{r}{0.994} & \multicolumn{1}{r}{} & \multicolumn{1}{r}{0.955} &
\multicolumn{1}{r}{0.996}\\
\multicolumn{1}{r}{} & \multicolumn{1}{r}{100} & \multicolumn{1}{r}{0.921} &
\multicolumn{1}{r}{0.995} & \multicolumn{1}{r}{} & \multicolumn{1}{r}{0.946} &
\multicolumn{1}{r}{0.995}\\
\multicolumn{1}{r}{} & \multicolumn{1}{r}{200} & \multicolumn{1}{r}{0.914} &
\multicolumn{1}{r}{0.996} & \multicolumn{1}{r}{} & \multicolumn{1}{r}{0.945} &
\multicolumn{1}{r}{0.998}\\\hline
\end{tabular}
\caption{Efficiencies of estimators}\label{tabEfiMulti}%
\end{table}%

\end{center}

It is seen that DCML is able to substantially increase the efficiency of S-E,
especially for $p=2.$ The efficiency for location is much higher than for
scatter. The fact that the efficiency of S-E increases with $p$ is also clear.
Actually, for $p=15$ the efficiency of S-E is $\geq0.96.$

\begin{figure}
[tbh]
\begin{center}
\includegraphics[
height=9.1116cm,
width=12.1166cm
]%
{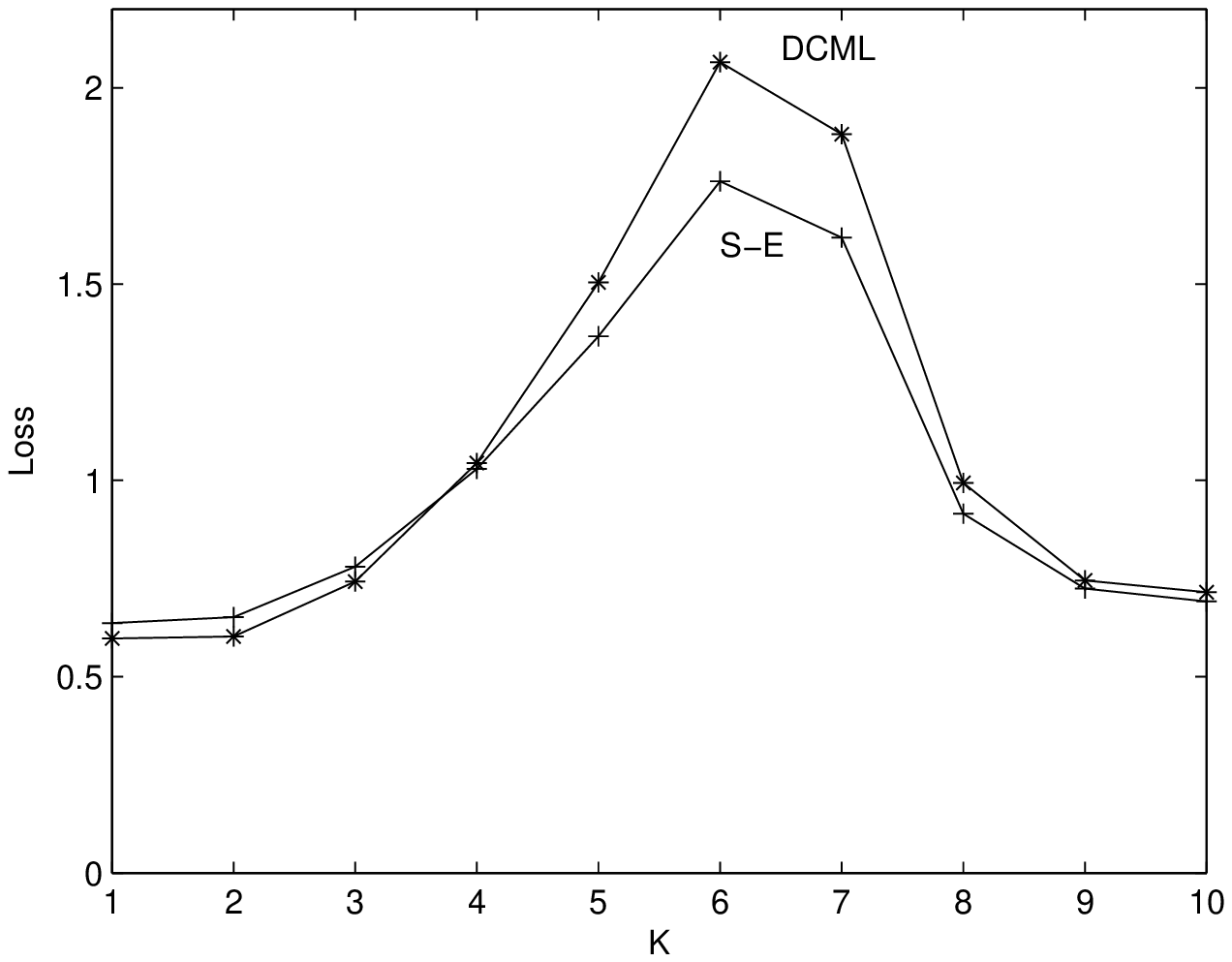}%
\caption{Losses of scatter matrices for $p=10,$ $n=100$ and 10\%
contamination, as a function of outlier size.}%
\label{figVMult_10}%
\end{center}
\end{figure}

Table \ref{tabContaMulti} shows the maximum mean losses for contamination rate
$\varepsilon=0.1.$ It is seen that in general the price for the increase in
efficiency is at worst a small increase of the maximum loss and at best a
decrease thereof. Figure \ref{figVMult_10} compares the losses of S-E and DCML
as a function of the outlier size $K$ for $\varepsilon=0.1.$%

\begin{center}%
\begin{table}[tbp] \centering
\begin{tabular}
[c]{cccccccc}\hline
&  & \textbf{\ } & \multicolumn{2}{c}{\textbf{\ }$\mathbf{\Sigma}$} &  &
\multicolumn{2}{c}{$\mathbf{\mu}$}\\\hline
$\varepsilon$ & $p$ & $n$ & S-E & DCML &  & S-E & DCML\\\cline{1-5}\cline{7-8}%
0.1 & \multicolumn{1}{r}{2} & \multicolumn{1}{r}{10} &
\multicolumn{1}{r}{0.91} & \multicolumn{1}{r}{1.03} & \multicolumn{1}{r}{} &
\multicolumn{1}{r}{0.32} & \multicolumn{1}{r}{0.34}\\
& \multicolumn{1}{r}{} & \multicolumn{1}{r}{20} & \multicolumn{1}{r}{0.51} &
\multicolumn{1}{r}{0.53} & \multicolumn{1}{r}{} & \multicolumn{1}{r}{0.18} &
\multicolumn{1}{r}{0.20}\\
& \multicolumn{1}{r}{} & \multicolumn{1}{r}{40} & \multicolumn{1}{r}{0.27} &
\multicolumn{1}{r}{0.31} & \multicolumn{1}{r}{} & \multicolumn{1}{r}{0.10} &
\multicolumn{1}{r}{0.11}\\
& \multicolumn{1}{r}{5} & \multicolumn{1}{r}{25} & \multicolumn{1}{r}{1.01} &
\multicolumn{1}{r}{1.03} & \multicolumn{1}{r}{} & \multicolumn{1}{r}{0.26} &
\multicolumn{1}{r}{0.27}\\
& \multicolumn{1}{r}{} & \multicolumn{1}{r}{50} & \multicolumn{1}{r}{0.65} &
\multicolumn{1}{r}{0.72} & \multicolumn{1}{r}{} & \multicolumn{1}{r}{0.17} &
\multicolumn{1}{r}{0.20}\\
& \multicolumn{1}{r}{} & \multicolumn{1}{r}{100} & \multicolumn{1}{r}{0.39} &
\multicolumn{1}{r}{0.48} & \multicolumn{1}{r}{} & \multicolumn{1}{r}{0.11} &
\multicolumn{1}{r}{0.13}\\
& \multicolumn{1}{r}{10} & \multicolumn{1}{r}{50} & \multicolumn{1}{r}{2.90} &
\multicolumn{1}{r}{3.26} & \multicolumn{1}{r}{} & \multicolumn{1}{r}{0.44} &
\multicolumn{1}{r}{0.51}\\
& \multicolumn{1}{r}{} & \multicolumn{1}{r}{100} & \multicolumn{1}{r}{1.82} &
\multicolumn{1}{r}{2.06} & \multicolumn{1}{r}{} & \multicolumn{1}{r}{0.28} &
\multicolumn{1}{r}{0.31}\\
& \multicolumn{1}{r}{} & \multicolumn{1}{r}{200} & \multicolumn{1}{r}{1.39} &
\multicolumn{1}{r}{1.71} & \multicolumn{1}{r}{} & \multicolumn{1}{r}{0.21} &
\multicolumn{1}{r}{0.27}\\\hline
\multicolumn{1}{r}{0.2} & \multicolumn{1}{r}{2} & \multicolumn{1}{r}{10} &
\multicolumn{1}{r}{1.42} & \multicolumn{1}{r}{1.49} & \multicolumn{1}{r}{} &
\multicolumn{1}{r}{0.50} & \multicolumn{1}{r}{0.51}\\
\multicolumn{1}{r}{} & \multicolumn{1}{r}{} & \multicolumn{1}{r}{20} &
\multicolumn{1}{r}{0.95} & \multicolumn{1}{r}{0.77} & \multicolumn{1}{r}{} &
\multicolumn{1}{r}{0.34} & \multicolumn{1}{r}{0.37}\\
\multicolumn{1}{r}{} & \multicolumn{1}{r}{} & \multicolumn{1}{r}{40} &
\multicolumn{1}{r}{0.68} & \multicolumn{1}{r}{0.50} & \multicolumn{1}{r}{} &
\multicolumn{1}{r}{0.28} & \multicolumn{1}{r}{0.27}\\
\multicolumn{1}{r}{} & \multicolumn{1}{r}{5} & \multicolumn{1}{r}{25} &
\multicolumn{1}{r}{3.60} & \multicolumn{1}{r}{3.43} & \multicolumn{1}{r}{} &
\multicolumn{1}{r}{0.97} & \multicolumn{1}{r}{1.19}\\
\multicolumn{1}{r}{} & \multicolumn{1}{r}{} & \multicolumn{1}{r}{50} &
\multicolumn{1}{r}{2.49} & \multicolumn{1}{r}{3.39} & \multicolumn{1}{r}{} &
\multicolumn{1}{r}{0.67} & \multicolumn{1}{r}{0.87}\\
\multicolumn{1}{r}{} & \multicolumn{1}{r}{} & \multicolumn{1}{r}{100} &
\multicolumn{1}{r}{2.20} & \multicolumn{1}{r}{2.13} & \multicolumn{1}{r}{} &
\multicolumn{1}{r}{0.56} & \multicolumn{1}{r}{0.73}\\
\multicolumn{1}{r}{} & \multicolumn{1}{r}{10} & \multicolumn{1}{r}{50} &
\multicolumn{1}{r}{11.21} & \multicolumn{1}{r}{11.15} & \multicolumn{1}{r}{} &
\multicolumn{1}{r}{2.49} & \multicolumn{1}{r}{2.88}\\
\multicolumn{1}{r}{} & \multicolumn{1}{r}{} & \multicolumn{1}{r}{100} &
\multicolumn{1}{r}{6.46} & \multicolumn{1}{r}{6.50} & \multicolumn{1}{r}{} &
\multicolumn{1}{r}{1.65} & \multicolumn{1}{r}{1.95}\\
\multicolumn{1}{r}{} & \multicolumn{1}{r}{} & \multicolumn{1}{r}{200} &
\multicolumn{1}{r}{5.71} & \multicolumn{1}{r}{5.72} & \multicolumn{1}{r}{} &
\multicolumn{1}{r}{1.52} & \multicolumn{1}{r}{1.78}\\\hline
\end{tabular}
\caption{Simulation: maximum mean losses of estimators}\label{tabContaMulti}%
\end{table}%

\end{center}

\subsection{Breakdown point}

It is easy to show that the replacement breakdown point of the DCML estimators
is that of the initial ones. We give the details for $\widehat{\mathbf{\Sigma
}},$ the case of $\widehat{\mathbf{\mu}}$ being similar. Consider a data set
$\mathbf{X=\{x}_{i}$ $i=1,...,n\}.$ Let $m$ be such that $\varepsilon=m/n$ is
less than the breakdown point $\varepsilon^{\ast}$ of the initial estimator
$\widehat{\mathbf{\Sigma}}_{0}.$ Let $\mathbf{X}^{\ast}$ be a data set that
coincides with $\mathbf{X}$ except for $m$ elements which are arbitrary. We
have to prove that, as a function of $\mathbf{X}^{\ast},$ the largest
eigenvalue $\lambda_{\max}$ of $\widehat{\mathbf{\Sigma}}$ is bounded, and the
smallest one $\lambda_{\min}$ is bounded away from zero. We know that this
property holds for $\widehat{\mathbf{\Sigma}}_{0}.$ Since by
(\ref{defDIstSIgma})
\[
\log|\widehat{\mathbf{\Sigma}}\mathbf{|-}\log|\widehat{\mathbf{\Sigma}}%
_{0}\mathbf{|+}\mathrm{trace}\left(  \widehat{\mathbf{\Sigma}}^{-1}%
\widehat{\mathbf{\Sigma}}_{0}\right)  -p\leq\delta,
\]
it follows from the \textquotedblleft trace\textquotedblright\ term that
$\lambda_{\min}$ cannot tend to zero, and then it follows from the
\textquotedblleft log\textquotedblright\ term that $\lambda_{\max}$ cannot
tend to infinity. $\blacksquare$

\section{\bigskip Real data\label{SecReal}}

In this section we apply the estimators to two published data sets.

\subsection{Regression}

We consider the well-known stack loss data set with $n=21$ and $p=3$ plus
intercept. Lacking a \textquotedblleft true model\textquotedblright\ we have
to employ alternative criteria for robustness and efficiency.

There seems to be a general agreement to consider observations 1, 3, 4 and 21
as atypical; see (Rousseew and Leroy, 1987). Call \textquotedblleft good
data\textquotedblright\ the data set without \{1,3,4,21\}. The estimators were
first computed using the good data, and the root mean squared prediction
errors (RMSE: square root of the mean of the squared residuals) was computed
for the same data. The comparison with LS was employed as a surrogate
criterion for efficiency. For a surrogate criterion for robustness, the
estimators were then computed for the whole data set, and the RMSE again
computed \emph{only} for the good data. Table \ref{tabStack} shows the results.

\begin{center}
\bigskip%
\begin{table}[tbp] \centering
\begin{tabular}
[c]{lcccccc}\hline
Computed with & LS & S-E & MM & G-Y & B-S & DCML\\\hline
\multicolumn{1}{r}{Good data} & \multicolumn{1}{r}{1.095} &
\multicolumn{1}{r}{1.416} & \multicolumn{1}{r}{1.126} & 1.095 & 1.095 &
\multicolumn{1}{r}{1.095}\\
\multicolumn{1}{r}{Whole data} & \multicolumn{1}{r}{1.921} &
\multicolumn{1}{r}{1.143} & \multicolumn{1}{r}{1.100} & 1.322 & 1.484 &
\multicolumn{1}{r}{1.164}\\\hline
\end{tabular}
\caption{Stack los data: prediction RMSEs of estimators for ``good''  data}\label{tabStack}%
\end{table}%

\end{center}

The first row shows that G-Y, B-S and DCML are here \textquotedblleft fully
efficient\textquotedblright, S-E is rather inefficient, and MM has a high
efficiency. The second row shows S-E, MM and DCML as most robust, followed by
G-Y, and B-S as the less robust one.

The behavior of S-E is puzzling. It gives zero weights to some
\textquotedblleft good\textquotedblright\ observations. The estimator was
recomputed several times to rule out the effect of the subsampling.

\subsection{Multivariate estimation}

Here we choose the Philips Mecoma data, employed in Problem 1 in (Rousseeuw
and Van Driessen, 1999), with $n=677$ and $p=9.$ Plotting the Mahalanobis
distances from the S estimator shows a number of clear outliers, the sequence
with indexes between 491 and 565 being the most outstanding ones. We defined
as \textquotedblleft bad data\textquotedblright\ the observations with
Mahalanobis distances larger than 60, which yielded 80 observations. Lacking a
criterion similar to prediction error like in the former example, we defined
as the \textquotedblleft true parameters\textquotedblright\ the MLE (mean and
covariance matrix) applied to the \textquotedblleft good\textquotedblright%
\ data, which will be called $\mathbf{\mu}_{\mathrm{good}}$ and
$\mathbf{\Sigma}_{\mathrm{good}},$ respectively.

We then computed, as above, the estimators based on the \textquotedblleft
good\textquotedblright\ data and their Kullback-Leibler distances to the
\textquotedblleft truth\textquotedblright; and then did the same for the
estimators based on the whole data. Namely, we computed
\[
d=\mathrm{trace}\left(  \mathbf{\Sigma}_{\mathrm{god}}^{-1}\mathbf{V}\right)
-p-\log|\mathbf{\Sigma}_{\mathrm{good}}^{-1}\mathbf{V}|
\]
for each scatter estimator $\mathbf{V}$, and
\[
d=\left(  \mathbf{t-\mu}_{\mathrm{good}}\right)  ^{\prime}\mathbf{C}%
_{\mathrm{good}}^{-1}\left(  \mathbf{t-\mu}_{\mathrm{good}}\right)
\]
for each location estimator $\mathbf{t.}$ Table \ref{tabPhilips} shows the results.

\begin{center}%
\begin{table}[tbp] \centering
\begin{tabular}
[c]{ccccc}\hline
\textbf{\ } & \textbf{\ }Computed with & MLE & S-E & DCML\\\hline
Scatter & Good data & \multicolumn{1}{r}{} & \multicolumn{1}{r}{0.381} &
\multicolumn{1}{r}{0.286}\\
& Whole data & \multicolumn{1}{r}{6.282} & \multicolumn{1}{r}{0.381} &
\multicolumn{1}{r}{0.322}\\\hline
Location & Good data & \multicolumn{1}{r}{} & \multicolumn{1}{r}{0.051} &
\multicolumn{1}{r}{0.039}\\
& Whole data & \multicolumn{1}{r}{1.067} & \multicolumn{1}{r}{0.051} &
\multicolumn{1}{r}{0.044}\\\hline
\end{tabular}
\caption{Philips data: Kullback-Leibler distances between estimators and
``true values''}\label{tabPhilips}%
\end{table}%

\end{center}

It is seen that here DCML\ outperforms S-E in all cases.

\medskip

\textbf{References}

Bondell, H.D. and Stefanski, L.A. (2013). \textquotedblleft Efficient Robust
Regression via Two-Stage Generalized Empirical Likelihood\textquotedblright,
\emph{Journal of the American Statistical Association,} 108, 644-655

Davies, P.L.\ (1987), \textquotedblleft Asymptotic Behavior of S-estimators of
Multivariate Location Parameters and Dispersion Matrices\textquotedblright,
\textit{The} \textit{Annals of Statistics}, 15, 1269--1292.

Gervini, D. and Yohai, V.J. (2002), \textquotedblleft A Class of Robust and
Fully Efficient Regression Estimators\textquotedblright, \textit{The Annals of
Statistics}, 30, 583-616.

Hampel, F.R. (1971), \textquotedblleft A General Definition of Qualitative
Robustness\textquotedblright, \textit{The Annals of Mathematical
Statistics\/,} 42, 1887--1896.

Huber, P.J. (1964), \textquotedblleft Robust Estimation of a Location
Parameter\textquotedblright, \textit{The Annals of Mathematical Statistics,}
35\textbf{,} 73-101.

Maronna, R.A., Martin, R.D. and Yohai, V.J.(2006). \textit{Robust Statistics:
Theory and Methods}, New York: John Wiley and Sons.

Rocke, D.M. (1996), \textquotedblleft Robustness Properties of $S$-estimators
of Multivariate Location and Shape in High Dimension\textquotedblright,
\textit{The Annals of Statistics, }24\textbf{, }1327-1345.

Rousseeuw, P.J and Leroy, A.M. (1987), \emph{Robust Regression and Outlier
Detection}, New York: John Wiley and Sons.

Rousseeuw, P.J. and Van Driessen, K. (1999), \textquotedblleft A Fast
Algorithm for the Minimum Covariance Determinant Estimator\textquotedblright,
\emph{Technometrics} 41, 212-223.

Rousseeuw, P.J and Yohai, V.J. (1984), \textquotedblleft Robust Regression by
Means of S--estimators\textquotedblright, In \textit{Robust and Nonlinear Time
Series\/}, J. Franke, W. H\"{a}rdle y R. D. Martin (eds.). Lectures Notes in
Statistics, 26, 256--272, New York: Springer Verlag.

Yohai, V.J. (1987), \textquotedblleft High Breakdown--point and High
Efficiency Estimates for Regression\textquotedblright, \textit{The Annals of
Statistics\/,} 15, 642--65.

Yohai, V.J., Stahel, W.A. and Zamar, R.H. (1991), \textquotedblleft A
Procedure for Robust Estimation and Inference in Linear
Regression\textquotedblright, In \textit{Directions in Robust Statistics and
Diagnostics (Part II), }W. Stahel and S. Weisberg, eds., The IMA Volumes in
Mathematics and is Applications, New York: Springer Verlag, 365-374.

Yohai, V.J. and Zamar, R.H. (1988), \textquotedblleft High Breakdown Estimates
of Regression by Means of the Minimization of an Efficient
Scale\textquotedblright, \textit{Journal of the American Statistical
Association\/,} 83, 406--413.

\end{document}